# Mannheim Partner $D$-Curves in Euclidean 3-space $E^3$

**Mustafa Kazaz[a], H. Hüseyin Uğurlu[b], Mehmet Önder[a], Tanju Kahraman[a]**
[a] Celal Bayar University, Department of Mathematics, Faculty of Arts and Sciences, , Manisa, Turkey.
E-mails: mustafa.kazaz@bayar.edu.tr, mehmet.onder@bayar.edu.tr
[b] Gazi University, Gazi Faculty of Education, Department of Secondary Education Science and Mathematics Teaching, Mathematics Teaching Program, Ankara, Turkey. E-mail: hugurlu@gazi.edu.tr

**Abstract**
In this paper we consider the idea of Mannheim partner curves for curves lying on surfaces and by considering the Darboux frames of them we define these curves as Mannheim partner $D$-curves and give the characterizations for these curves. We also find the relations between the geodesic curvatures, the normal curvatures and the geodesic torsions of these associated curves. Furthermore, we show that the definition and the characterizations of Mannheim partner $D$-curves include those of Mannheim partner curves in some special cases.

## 1. Introduction

Associated curves, the curves for which at the corresponding points of the curves one of the Frenet vectors of a curve coincides with the one of the Frenet vectors of the other curves, are very interesting and an important problem of the fundamental theory and the characterizations of space curves. The well-known of such curves is Bertrand curve which is characterized as a kind of corresponding relation between the two curves. The relation is that the principal normal of a curve is the principal normal of another curve i.e, the Bertrand curve is a curve which shares the normal line with another curve. Over years many mathematicians have studied on Bertrand curves in different spaces and consider the properties of these curves[1,2,3,5,6]. By considering the frame of the ruled surface, Ravani and Ku extended the notion of Bertrand curve to the ruled surfaces and named as Bertrand offsets[14]. The corresponding characterizations of the Bertrand offsets of timelike ruled surface were given by Kurnaz[9].

Recently, a new definition of the associated curves was given by Liu and Wang[10,16]. They called these new curves as Mannheim partner curves: Let $x$ and $x_1$ be two curves in the three dimensional Euclidean $E^3$. If there exists a corresponding relationship between the space curves $x$ and $x_1$ such that, at the corresponding points of the curves, the principal normal lines of $x$ coincides with the binormal lines of $x_1$, then $x$ is called a Mannheim curve, and $x_1$ is called a Mannheim partner curve of $x$. The pair $\{x, x_1\}$ is said to be a Mannheim pair. They showed that the curve $x_1(s_1)$ is the Mannheim partner curve of the curve $x(s)$ if and only if the curvature $\kappa_1$ and the torsion $\tau_1$ of $x_1(s_1)$ satisfy the following equation

$$\dot{\tau} = \frac{d\tau}{ds_1} = \frac{\kappa_1}{\lambda}(1 + \lambda^2 \tau_1^2)$$

for some non-zero constant $\lambda$. They also study the Mannheim curves in Minkowski 3-space[10,16]. Similar to the Bertrand offsets, Orbay and others have defined and characterized the Mannheim offsets of ruled surfaces[13]. The corresponding characterizations of the Mannheim offsets of timelike and spacelike ruled surfaces in Minkowski 3-space have given by Kazaz, Ugurlu and Onder[7,8].

In this paper we consider the notion of the Mannheim partner curve for the curves lying on the surfaces. We call these new associated curves as Mannheim partner $D$-curves and by using the Darboux frame of the curves we give the definition and the characterizations of these curves.

## 2. Darboux Frame of a Curve Lying on a Surface

Let $S$ be an oriented surface in three-dimensional Euclidean space $E^3$ and let consider a curve $x(s)$ lying on $S$ fully. Since the curve $x(s)$ is also in space, there exists Frenet frame $\{\vec{T}, \vec{N}, \vec{B}\}$ at each points of the curve where $\vec{T}$ is unit tangent vector, $\vec{N}$ is principal normal vector and $\vec{B}$ is binormal vector, respectively. The Frenet equations of the curve $x(s)$ is given by

$$\vec{T}' = \kappa \vec{N}$$
$$\vec{N}' = -\kappa \vec{T} + \tau \vec{B}$$
$$\vec{B}' = -\tau \vec{N}$$

where $\kappa$ and $\tau$ are curvature and torsion of the curve $x(s)$, respectively.

Since the curve $x(s)$ lies on the surface $S$ there exists another frame of the curve $x(s)$ which is called Darboux frame and denoted by $\{\vec{T}, \vec{g}, \vec{n}\}$. In this frame $\vec{T}$ is the unit tangent of the curve, $\vec{n}$ is the unit normal of the surface $S$ and $\vec{g}$ is a unit vector given by $\vec{g} = \vec{n} \times \vec{T}$. Since the unit tangent $\vec{T}$ is common in both Frenet frame and Darboux frame, the vectors $\vec{N}$, $\vec{B}$, $\vec{g}$ and $\vec{n}$ lie on the same plane. So that the relations between these frames can be given as follows

$$\begin{bmatrix} \vec{T} \\ \vec{g} \\ \vec{n} \end{bmatrix} = \begin{bmatrix} 1 & 0 & 0 \\ 0 & \cos\varphi & \sin\varphi \\ 0 & -\sin\varphi & \cos\varphi \end{bmatrix} \begin{bmatrix} \vec{T} \\ \vec{N} \\ \vec{B} \end{bmatrix},$$

where $\varphi$ is the angle between the vectors $\vec{g}$ and $\vec{N}$. The derivative formulae of the Darboux frame is

$$\begin{bmatrix} \dot{\vec{T}} \\ \dot{\vec{g}} \\ \dot{\vec{n}} \end{bmatrix} = \begin{bmatrix} 0 & k_g & k_n \\ -k_g & 0 & \tau_g \\ -k_n & -\tau_g & 0 \end{bmatrix} \begin{bmatrix} \vec{T} \\ \vec{g} \\ \vec{n} \end{bmatrix} \qquad (1)$$

where $k_g$, $k_n$ and $\tau_g$ are called the geodesic curvature, the normal curvature and the geodesic torsion, respectively. Here and in the following, we use "dot" to denote the derivative with respect to the arc length parameter of a curve.

The relations between geodesic curvature, normal curvature, geodesic torsion and $\kappa$, $\tau$ are given as follows

$$k_g = \kappa \cos\varphi, \quad k_n = \kappa \sin\varphi, \quad \tau_g = \tau + \frac{d\varphi}{ds}. \qquad (2)$$

Furthermore, the geodesic curvature $k_g$ and geodesic torsion $\tau_g$ of the curve $x(s)$ can be calculated as follows

$$k_g = \left\langle \frac{d\vec{x}}{ds}, \frac{d^2\vec{x}}{ds^2} \times \vec{n} \right\rangle, \quad \tau_g = \left\langle \frac{d\vec{x}}{ds}, \vec{n} \times \frac{d\vec{n}}{ds} \right\rangle \qquad (3)$$

In the differential geometry of surfaces, for a curve $x(s)$ lying on a surface $S$ the followings are well-known

    **i)** $x(s)$ is a geodesic curve $\Leftrightarrow k_g = 0$,

    **ii)** $x(s)$ is an asymptotic line $\Leftrightarrow k_n = 0$,

    **iii)** $x(s)$ is a principal line $\Leftrightarrow \tau_g = 0$ [12].

Through every point of the surface passes a geodesic in every direction. A geodesic is uniquely determined by an initial point and tangent at that point. All straight lines on a surface are geodesics. Along all curved geodesics the principal normal coincides with the surface normal. Along asymptotic lines osculating planes and tangent planes coincide, along geodesics they are normal. Through a point of a nondevelopable surface pass two asymptotic lines which can be real or imaginary [15].

## 3. Mannheim Partner $D$-Curves in Euclidean 3-space $E^3$

In this section, by considering the Darboux frame, we define Mannheim partner $D$-curves and give the characterizations of these curves.

**Definition 1.** Let $S$ and $S_1$ be oriented surfaces in three-dimensional Euclidean space $E^3$ and let consider the arc-length parameter curves $x(s)$ and $x_1(s_1)$ lying fully on $S$ and $S_1$, respectively. Denote the Darboux frames of $x(s)$ and $x_1(s_1)$ by $\{\vec{T}, \vec{g}, \vec{n}\}$ and $\{\vec{T_1}, \vec{g_1}, \vec{n_1}\}$, respectively. If there exists a corresponding relationship between the curves $x$ and $x_1$ such that, at the corresponding points of the curves, the Darboux frame element $\vec{g}$ of $x$ coincides with the Darboux frame element $\vec{n_1}$ of $x_1$, then $x$ is called a Mannheim $D$-curve, and $x_1$ is a Mannheim partner $D$-curve of $x$. Then, the pair $\{x, x_1\}$ is said to be a Mannheim $D$-pair. If there exist such curves lying on the oriented surfaces $S$ and $S_1$, respectively, we call the pair $\{S, S_1\}$ as Mannheim pair surfaces.

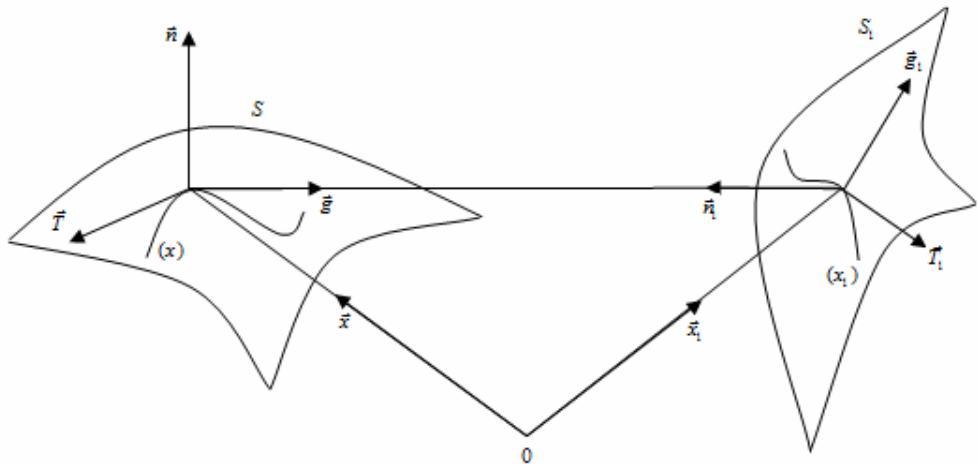

**Fig. 1** Mannheim partner $D$-curves

**Theorem 1.** *Let $S$ be an oriented surface and $x(s)$ be a Mannheim $D$-curve in $E^3$ with arc length parameter $s$ fully lying on $S$. If $S_1$ is another oriented surface and $x_1(s_1)$ is a curve with arc length parameter $s_1$ fully lying on $S_1$, then $x_1(s_1)$ is Mannheim partner $D$-curve of $x(s)$ if and only if the normal curvature $k_n$ of $x(s)$ and the geodesic curvature $k_{g_1}$, the normal curvature $k_{n_1}$ and the geodesic torsion $\tau_{g_1}$ of $x_1(s_1)$ satisfy the following equation*

$$-\lambda \dot{\tau}_{g_1} = \left( \frac{(1-\lambda k_{n_1})^2 + \lambda^2 \tau_{g_1}^2}{(1-\lambda k_{n_1})} \right) \left( \frac{\lambda k_{n_1} - 1}{\cos\theta} k_n + k_{g_1} \right) + \frac{\lambda^2 \tau_{g_1} \dot{k}_{n_1}}{1 - \lambda k_{n_1}}$$

*for some nonzero constants $\lambda$, where $\theta$ is the angle between the tangent vectors $\vec{T}$ and $\vec{T}_1$ at the corresponding points of $x$ and $x_1$.*

**Proof:** Suppose that $S$ is an oriented surface and $x(s)$ is a Mannheim $D$-curve fully lying on $S$. Denote the Darboux frames of $x(s)$ and $x_1(s_1)$ by $\{\vec{T}, \vec{g}, \vec{n}\}$ and $\{\vec{T}_1, \vec{g}_1, \vec{n}_1\}$, respectively. Then by the definition we can assume that

$$\vec{x}(s_1) = \vec{x}_1(s_1) + \lambda(s_1)\vec{n}_1(s_1), \tag{4}$$

for some function $\lambda(s_1)$. By taking derivative of (4) with respect to $s_1$ and applying the Darboux formulas (1) we have

$$\vec{T}\frac{ds}{ds_1} = (1-\lambda k_{n_1})\vec{T}_1 + \dot{\lambda}\vec{n}_1 - \lambda \tau_{g_1} \vec{g}_1. \tag{5}$$

Since the direction of $\vec{n}_1$ coincides with the direction of $\vec{g}$, we get

$$\dot{\lambda}(s_1) = 0.$$

This means that $\lambda$ is a nonzero constant. Thus, the equality (5) can be written as follows

$$\vec{T}\frac{ds}{ds_1} = (1-\lambda k_{n_1})\vec{T}_1 - \lambda \tau_{g_1} \vec{g}_1. \tag{6}$$

On the other hand we have

$$\vec{T} = \cos\theta \vec{T}_1 + \sin\theta \vec{g}_1, \tag{7}$$

where $\theta$ is the angle between the tangent vectors $\vec{T}$ and $\vec{T}_1$ at the corresponding points of $x$ and $x_1$. By differentiating this last equation with respect to $s_1$, we get

$$(k_g \vec{g} + k_n \vec{n})\frac{ds}{ds_1} = -(\dot{\theta} + k_{g_1})\sin\theta \vec{T}_1 + (\dot{\theta} + k_{g_1})\cos\theta \vec{g}_1 + (k_{n_1}\cos\theta + \tau_{g_1}\sin\theta)\vec{n}_1. \tag{8}$$

From this equation and the fact that

$$\vec{n} = \sin\theta \vec{T}_1 - \cos\theta \vec{g}_1, \tag{9}$$

we get

$$(k_n \sin\theta \vec{T}_1 - k_n \cos\theta \vec{g}_1 + k_g \vec{g})\frac{ds}{ds_1} = -(\dot{\theta} + k_{g_1})\sin\theta \vec{T}_1 + (\dot{\theta} + k_{g_1})\cos\theta \vec{g}_1 \\ + (k_{n_1}\cos\theta + \tau_{g_1}\sin\theta)\vec{n}_1 \tag{10}$$

Since the direction of $\vec{n}_1$ is coincident with $\vec{g}$ we have

$$\dot{\theta} = -\left( k_n \frac{ds}{ds_1} + k_{g_1} \right). \tag{11}$$

From (6) and (7) and notice that $\vec{T}_1$ is orthogonal to $\vec{g}_1$ we obtain

$$\frac{ds}{ds_1} = \frac{1-\lambda k_{n_1}}{\cos\theta} = -\frac{\lambda \tau_{g_1}}{\sin\theta}. \tag{12}$$

Equality (12) gives us

$$\tan\theta = -\frac{\lambda \tau_{g_1}}{1-\lambda k_{n_1}}. \tag{13}$$

By taking the derivative of this equation and applying (11) we get

$$-\lambda \dot{\tau}_{g_1} = \left(\frac{(1-\lambda k_{n_1})^2 + \lambda^2 \tau_{g_1}^2}{(1-\lambda k_{n_1})}\right)\left(\frac{\lambda k_{n_1}-1}{\cos\theta}k_n + k_{g_1}\right) + \frac{\lambda^2 \tau_{g_1} \dot{k}_{n_1}}{1-\lambda k_{n_1}}, \tag{14}$$

that is desired.

Conversely, assume that the equation (14) holds for some nonzero constants $\lambda$. Then by using (12) and (13), (14) gives us

$$-k_n\left(\frac{ds}{ds_1}\right)^3 = -\lambda \dot{\tau}_{g_1}(1-\lambda k_{n_1}) - \lambda^2 \tau_{g_1} \dot{k}_{n_1} + \left((1-\lambda k_{n_1})^2 + \lambda^2 \tau_{g_1}^2\right)k_{g_1} \tag{15}$$

Let define a curve

$$\vec{x}(s_1) = \vec{x}_1(s_1) + \lambda \vec{n}_1(s_1). \tag{16}$$

where $\lambda$ is a non-zero constant. We will prove that $x$ is a Mannheim $D$-curve and $x_1$ is the Mannheim partner $D$-curve of $x$. By taking the derivative of (16) with respect to $s_1$ twice, we get

$$\vec{T}\frac{ds}{ds_1} = (1-\lambda k_{n_1})\vec{T}_1 - \lambda \tau_{g_1}\vec{g}_1, \tag{17}$$

and

$$(k_g \vec{g} + k_n \vec{n})\left(\frac{ds}{ds_1}\right)^2 + \vec{T}\frac{d^2 s}{ds_1^2} = (-\lambda \dot{k}_{n_1} + \lambda \tau_{g_1} k_{g_1})\vec{T}_1 + \left((1-\lambda k_{n_1})k_{g_1} - \lambda \dot{\tau}_{g_1}\right)\vec{g}_1 \\ + \left((1-\lambda k_{n_1})k_{n_1} - \lambda \tau_{g_1}^2\right)\vec{n}_1 \tag{18}$$

respectively. Taking the cross product of (17) with (18) we have

$$\left[k_g \vec{n} - k_n \vec{g}\right]\left(\frac{ds}{ds_1}\right)^3 = \left(-\lambda \tau_{g_1} k_{n_1}(1-\lambda k_{n_1}) + \lambda^2 \tau_{g_1}^3\right)\vec{T}_1 - \left((1-\lambda k_{n_1})^2 k_{n_1} - \lambda \tau_{g_1}^2(1-\lambda k_{n_1})\right)\vec{g}_1 \\ + \left(k_{g_1}(1-\lambda k_{n_1})^2 - \lambda \dot{\tau}_{g_1}(1-\lambda k_{n_1}) - \lambda^2 \tau_{g_1} \dot{k}_{n_1} + \lambda^2 \tau_{g_1}^2 k_{g_1}\right)\vec{n}_1 \tag{19}$$

By substituting (15) in (19) we get

$$\left[k_g \vec{n} - k_n \vec{g}\right]\left(\frac{ds}{ds_1}\right)^3 = \left(-\lambda \tau_{g_1} k_{n_1}(1-\lambda k_{n_1}) + \lambda^2 \tau_{g_1}^3\right)\vec{T}_1 \\ - \left(k_{n_1}(1-\lambda k_{n_1})^2 - \lambda \tau_{g_1}^2(1-\lambda k_{n_1})\right)\vec{g}_1 \\ - k_n\left(\frac{ds}{ds_1}\right)^3 \vec{n}_1 \tag{20}$$

Taking the cross product of (17) with (20) we have

$$-\left[k_g \vec{g} + k_n \vec{n}\right]\left(\frac{ds}{ds_1}\right)^4 = k_n\left(\frac{ds}{ds_1}\right)^3 \lambda \tau_{g_1}\vec{T}_1 + k_n\left(\frac{ds}{ds_1}\right)^3(1-\lambda k_{n_1})\vec{g}_1 \\ + \left((1-\lambda k_{n_1})^2 + \lambda^2 \tau_{g_1}^2\right)\left(\lambda \tau_{g_1}^2 - k_{n_1}(1-\lambda k_{n_1})\right)\vec{n}_1 \tag{21}$$

From (20) and (21) we have

$$\left(k_g^2 + k_n^2\right)\left(\frac{ds}{ds_1}\right)^4 \vec{n} = \left[k_g \frac{ds}{ds_1}\left(-\lambda\tau_{g_1} k_{n_1}(1-\lambda k_{n_1}) + \lambda^2\tau_{g_1}^3\right) - \lambda\tau_{g_1} k_n^2 \left(\frac{ds}{ds_1}\right)^3\right]\vec{T_1}$$

$$- \left[k_g \frac{ds}{ds_1}\left(k_{n_1}(1-\lambda k_{n_1})^2 - \lambda\tau_{g_1}^2(1-\lambda k_{n_1})\right) + (1-\lambda k_{n_1})k_n^2\left(\frac{ds}{ds_1}\right)^3\right]\vec{g_1} \quad (22)$$

$$- \left[k_n k_g \left(\frac{ds}{ds_1}\right)^4 + k_n\left((1-\lambda k_{n_1})^2 + \lambda^2\tau_{g_1}^2\right)\left(\lambda\tau_{g_1}^2 - k_{n_1}(1-\lambda k_{n_1})\right)\right]\vec{n_1}$$

Furthermore, from (17) and (20) we get

$$\begin{cases} \left(\dfrac{ds}{ds_1}\right)^2 = (1-\lambda k_{n_1})^2 + \lambda^2\tau_{g_1}^2, \\ k_g\left(\dfrac{ds}{ds_1}\right)^2 = k_{n_1}(1-\lambda k_{n_1}) - \lambda\tau_{g_1}^2, \end{cases} \quad (23)$$

respectively. Substituting (23) in (22) we obtain

$$\left(k_g^2 + k_n^2\right)\left(\frac{ds}{ds_1}\right)^4 \vec{n} = \left[k_g \frac{ds}{ds_1}\left(-\lambda\tau_{g_1} k_{n_1}(1-\lambda k_{n_1}) + \lambda^2\tau_{g_1}^3\right) - \lambda\tau_{g_1} k_n^2 \left(\frac{ds}{ds_1}\right)^3\right]\vec{T_1}$$

$$- \left[k_g \frac{ds}{ds_1}\left(k_{n_1}(1-\lambda k_{n_1})^2 - \lambda\tau_{g_1}^2(1-\lambda k_{n_1})\right) + (1-\lambda k_{n_1})k_n^2\left(\frac{ds}{ds_1}\right)^3\right]\vec{g_1} \quad (24)$$

Equality (17) and (24) shows that the vectors $\vec{T}$ and $\vec{n}$ lie on the plane $sp\{\vec{T_1}, \vec{g_1}\}$. So, at the corresponding points of the curves, the Darboux frame element $\vec{g}$ of $x$ coincides with the Darboux frame element $\vec{n_1}$ of $x_1$, i.e, the curves $x$ and $x_1$ are Mannheim $D$-pair curves.

Let now give the characterizations of Mannheim partner $D$-curves in some special cases. Assume that $x(s)$ be an asymptotic Mannheim $D$-curve. Then, from (14) we have the following special cases:

**i)** Consider that $x_1(s_1)$ is a geodesic curve. Then $x_1(s_1)$ is Mannheim partner $D$-curve of $x(s)$ if and only if the following equation holds,

$$\dot{\tau}_{g_1} = -\frac{\lambda\tau_{g_1}\dot{k}_{n_1}}{1-\lambda k_{n_1}}.$$

**ii)** Assume that $x_1(s_1)$ is also an asymptotic line. Then $x_1(s_1)$ is Mannheim partner $D$-curve of $x(s)$ if and only if the geodesic curvature $k_{g_1}$ and the geodesic torsion $\tau_{g_1}$ of $x_1(s_1)$ satisfy the following equation,

$$-\lambda\dot{\tau}_{g_1} = (1+\lambda^2\tau_{g_1}^2)k_{g_1}.$$

In this case, the Frenet frame of the curve $x_1(s_1)$ coincides with its Darboux frame. From (2) we have $k_{g_1} = \kappa_1$ and $\tau_{g_1} = \tau_1$. So, the Mannheim partner $D$-curves become the Mannheim partner curves, i.e., if both $x(s)$ and $x_1(s_1)$ are asymptotic lines then, the definition and the characterizations of the Mannheim partner $D$-curves involve those of the Mannheim partner curves in Euclidean 3-space.

**iii)** If $x_1(s_1)$ is a principal line then $x_1(s_1)$ is Mannheim partner $D$-curve of $x(s)$ if and only if the geodesic curvature $k_{g_1} = 0$ i.e, $x_1(s_1)$ is also a geodesic curve or $k_{n_1} = 1/\lambda = const$.

**Theorem 2.** Let the pair $\{x, x_1\}$ be a Mannheim $D$-pair. Then the relation between geodesic curvature $k_g$, geodesic torsion $\tau_g$ of $x(s)$ and the normal curvature $k_{n_1}$, the geodesic torsion $\tau_{g_1}$ of $x_1(s_1)$ is given as follows

$$k_g - k_{n_1} = \lambda(k_g k_{n_1} - \tau_g \tau_{g_1}).$$

**Proof:** Let $x(s)$ ba a Mannheim $D$-curve and $x_1(s_1)$ be a Mannheim partner $D$-curve of $x(s)$. Then from (16) we can write

$$\vec{x}_1(s_1) = \vec{x}(s_1) - \lambda \vec{n}_1(s_1) \tag{25}$$

for some constants $\lambda$. By differentiating (25) with respect to $s_1$ we have

$$\vec{T}_1 = (1 + \lambda k_g)\frac{ds}{ds_1}\vec{T} - \lambda \tau_g \frac{ds}{ds_1}\vec{n} \tag{26}$$

By the definition we have

$$\vec{T}_1 = \cos\theta \vec{T} + \sin\theta \vec{n} \tag{27}$$

From (26) and (27) we obtain

$$\cos\theta = (1 + \lambda k_g)\frac{ds}{ds_1}, \quad \sin\theta = -\lambda \tau_g \frac{ds}{ds_1}. \tag{28}$$

Using (12) and (28) it is easily seen that

$$k_g - k_{n_1} = \lambda(k_g k_{n_1} - \tau_g \tau_{g_1}).$$

From Theorem 2, we obtain the following special cases.
Let the pair $\{x, x_1\}$ be a Mannheim $D$-pair. Then,

    **i)** if $x_1$ is an asymptotic line, then
$$k_g = -\lambda \tau_g \tau_{g_1}$$
    **ii)** if $x_1$ is a principal line, then
$$k_g - k_{n_1} = \lambda k_g k_{n_1}$$
    **iii)** if $x$ is a geodesic curve, then
$$k_{n_1} = \lambda \tau_g \tau_{g_1}$$
    **iv)** if $x$ is a principal line then
$$k_g - k_{n_1} = \lambda k_g k_{n_1}$$

**Theorem 3.** Let $\{x, x_1\}$ be Mannheim $D$-pair. Then the following relations hold:

    **i)** $k_{g_1} = -\left(k_n \dfrac{ds}{ds_1} + \dfrac{d\theta}{ds_1}\right)$

    **ii)** $\tau_g \dfrac{ds}{ds_1} = -k_{n_1} \sin\theta + \tau_{g_1} \cos\theta$

    **iii)** $k_g \dfrac{ds}{ds_1} = k_{n_1} \cos\theta + \tau_{g_1} \sin\theta$

**iv)** $\tau_{g_1} = (k_g \sin\theta + \tau_g \cos\theta) \dfrac{ds}{ds_1}$

**Proof: i)** By differentiating the equation $\langle \vec{T}, \vec{T_1} \rangle = \cos\theta$ with respect to $s_1$ we have

$$\left\langle (k_g \vec{g} + k_n \vec{n}) \frac{ds}{ds_1}, \vec{T_1} \right\rangle + \left\langle \vec{T}, k_{g_1} \vec{g_1} + k_{n_1} \vec{n_1} \right\rangle = -\sin\theta \frac{d\theta}{ds_1}.$$

Using the fact that the direction of $n_1$ coincides with the direction of $g$ and

$$\begin{cases} \vec{T_1} = \cos\theta \vec{T} + \sin\theta \vec{n}, \\ \vec{g_1} = \sin\theta \vec{T} - \cos\theta \vec{n}, \end{cases} \qquad (29)$$

we easily get that

$$k_{g_1} = -\left( k_n \frac{ds}{ds_1} + \frac{d\theta}{ds_1} \right).$$

**ii)** By differentiating the equation $\langle \vec{n}, \vec{n_1} \rangle = 0$ with respect to $s_1$ we have

$$\left\langle (-k_n \vec{T} - \tau_g \vec{g}) \frac{ds}{ds_1}, \vec{n_1} \right\rangle + \left\langle \vec{n}, -k_{n_1} \vec{T_1} - \tau_{g_1} \vec{g_1} \right\rangle = 0.$$

By (29) we obtain

$$\tau_g \frac{ds}{ds_1} = -k_{n_1} \sin\theta + \tau_{g_1} \cos\theta.$$

**iii)** By differentiating the equation $\langle \vec{T}, \vec{n_1} \rangle = 0$ with respect to $s_1$ we get

$$\left\langle (k_g \vec{g} + k_n \vec{n}) \frac{ds}{ds_1}, \vec{n_1} \right\rangle + \left\langle \vec{T}, -k_{n_1} \vec{T_1} - \tau_{g_1} \vec{g_1} \right\rangle = 0.$$

From (29) it follows that

$$k_g \frac{ds}{ds_1} = k_{n_1} \cos\theta + \tau_{g_1} \sin\theta.$$

**iv)** By differentiating the equation $\langle \vec{g}, \vec{g_1} \rangle = 0$ with respect to $s_1$ we obtain

$$\left\langle (-k_g \vec{T} + \tau_g \vec{n}) \frac{ds}{ds_1}, \vec{g_1} \right\rangle + \left\langle \vec{g}, -k_{g_1} \vec{T_1} + \tau_{g_1} \vec{n_1} \right\rangle = 0.$$

By considering (29) we get

$$\tau_{g_1} = (k_g \sin\theta + \tau_g \cos\theta) \frac{ds}{ds_1}.$$

Let now $x$ be a Mannheim $D$-curve and $x_1$ be a Mannheim partner $D$-curve of $x$. From the first equation of (3) and by using the fact that $n_1$ is coincident with $g$ we have

$$k_{g_1} = \langle \dot{\vec{x}}_1, \ddot{\vec{x}}_1 \times \vec{n}_1 \rangle = \langle \dot{\vec{x}}_1, \ddot{\vec{x}}_1 \times \vec{g} \rangle$$
$$= \left(\frac{ds}{ds_1}\right)^3 \left( -k_{n_1}(1+\lambda k_g)^2 - \lambda^2 \tau_g^2 k_n \right) + \left(\frac{ds}{ds_1}\right)^2 \left( \lambda \dot{\tau}_g (1+\lambda k_g) - \lambda^2 \tau_g \dot{k}_g \right) \qquad (30)$$

Then the relations between the geodesic curvature $k_{g_1}$ of $x_1(s_1)$ and the geodesic curvature $k_g$, the normal curvature $k_n$ and the geodesic torsion $\tau_g$ of $x(s)$ are given as follows,

If $k_g = 0$ then from (30) the geodesic curvature $k_{g_1}$ of $x_1(s_1)$ is

$$k_{g_1} = -\left(\frac{ds}{ds_1}\right)^3 (1+\lambda^2\tau_g^2)k_n + \left(\frac{ds}{ds_1}\right)^2 \lambda\dot{\tau}_g. \tag{31}$$

If $k_n = 0$ then the relation between $k_g$, $\tau_g$ and $k_{g_1}$ is

$$k_{g_1} = \lambda\left(\frac{ds}{ds_1}\right)^2 \left(\dot{\tau}_g(1+\lambda k_g) - \lambda\tau_g\dot{k}_g\right). \tag{32}$$

If $\tau_g = 0$ then, for the geodesic curvature $k_{g_1}$, we have

$$k_{g_1} = -\left(\frac{ds}{ds_1}\right)^3 (1+\lambda k_g)^2 k_n. \tag{33}$$

From (31), (32) and (33) we give the following corollary.

**Corollary 1.** *Let $x$ be a Mannheim $D$-curve and $x_1$ be a Mannheim partner $D$-curve of $x$. Then the relations between the geodesic curvature $k_{g_1}$ of $x_1(s_1)$ and the geodesic curvature $k_g$, the normal curvature $k_n$ and the geodesic torsion $\tau_g$ of $x(s)$ are given as follows,*

  *i) If $x$ is a geodesic curve, then the geodesic curvature $k_{g_1}$ of $x_1(s_1)$ is*

$$k_{g_1} = -\left(\frac{ds}{ds_1}\right)^3 (1+\lambda^2\tau_g^2)k_n + \left(\frac{ds}{ds_1}\right)^2 \lambda\dot{\tau}_g.$$

  *ii) If $x$ is an asymptotic line, then the equation of $k_{g_1}$ is*

$$k_{g_1} = \lambda\left(\frac{ds}{ds_1}\right)^2 \left(\dot{\tau}_g(1+\lambda k_g) - \lambda\tau_g\dot{k}_g\right).$$

  *iii) If $x$ is a principal line, then the geodesic curvature $k_{g_1}$ of $x_1(s_1)$ is*

$$k_{g_1} = -\left(\frac{ds}{ds_1}\right)^3 (1+\lambda k_g)^2 k_n.$$

Similarly, From the second equation of (3) and by using the fact that $g$ is coincident with $n_1$, the relation between the geodesic torsion $\tau_{g_1}$ of $x_1(s_1)$ and the geodesic torsion $\tau_g$ of $x(s)$ is given by

$$\tau_{g_1} = \left(\frac{ds}{ds_1}\right)^2 \tau_g. \tag{34}$$

Furthermore, by using (12), from (34) we have

$$\tau_g \tau_{g_1} = \frac{\sin^2\theta}{\lambda^2}. \tag{35}$$

Then, from (34) and (35) we can give the following corollary.

**Corollary 2.** *Let $x$ be a Mannheim $D$-curve and $x_1$ be a Mannheim partner $D$-curve of $x$. Then the relation between the geodesic torsion $\tau_{g_1}$ of $x_1(s_1)$ and the geodesic torsion $\tau_g$ of $x(s)$ is given by one of the followings,*

$$\tau_{g_1} = \left(\frac{ds}{ds_1}\right)^2 \tau_g \quad \text{or} \quad \tau_g\tau_{g_1} = \frac{\sin^2\theta}{\lambda^2}$$

*and so, the Mannheim partner $D$-curve $x_1$ is a principal line when the Mannheim $D$-curve $x$ is a principal line.*

Similarly, from (12) and (32) we get

$$\frac{\tau_g}{\tau_{g_1}} = \frac{\cos^2\theta}{(1-\lambda k_{n_1})^2}$$

Then, if $x_1(s_1)$ is an asymptotic curve, i.e., $k_{n_1}=0$, we have

$$\tau_g = \cos^2\theta \tau_{g_1} \qquad (36)$$

From (36) we have the following corollary.

**Corollary 3.** *Let $x$ be a Mannheim $D$-curve and $x_1$ be a Mannheim partner $D$-curve of $x$. If $x_1(s_1)$ is an asymptotic curve then the relation between the geodesic torsion $\tau_g$ of $x(s)$ and the geodesic torsion $\tau_{g_1}$ of $x_1(s_1)$ is given as follows,*

$$\tau_g = \cos^2\theta \tau_{g_1}$$

*where $\theta$ is the angle between the tangent vectors $T$ and $T_1$ at the corresponding points of $x$ and $x_1$.*

## 4. Conclusions

In this paper, the definition and characterizations of Mannheim partner $D$-curves are given which is a new study of associated curves lying on surfaces. It is shown that the definition and the characterizations of Mannheim partner $D$-curves include those of Mannheim partner curves in some special cases. Furthermore, the relations between the geodesic curvature, the normal curvature and the geodesic torsion of these curves are given.